\documentstyle[11pt]{article}
\begin{document}

\rightline{ITP-97-9E \hspace{12mm}}
\rightline{q-alg/9709035 \hspace{3mm}}
\vspace{30mm}

\noindent
\centerline{\bf  REPRESENTATIONS OF THE $q$-DEFORMED ALGEBRA
U$'_{\bf q}$(so$_{\bf 2,2}$)}
\vspace{18mm}

\centerline{A. U. KLIMYK }

\vspace{12mm}
\centerline{
Institute for Theoretical Physics of Ukrainian National Academy of
Sciences,}
\centerline
{Metrologichna 14-B, 252143, Kiev, Ukraine}

\vspace{24mm}

\begin{abstract}

The main aim of this paper is to give classes of irreducible
infinite dimensional representations and of irreducible $*$-representations
of the $q$-deformed algebra $U'_q({\rm so}_{2,2})$ which is a real form
of the non-standard deformation $U'_q({\rm so}_4)$ of the universal
enveloping algebra $U({\rm so}(4,{\bf C}))$. These representations are
described by two complex parameters and are obtained by "analytical
continuation" of the irreducible finite dimensional representations of the
algebra $U'_q({\rm so}_4)$ in the basis corresponding to the reduction from
$U'_q({\rm so}_4)$ to $U({\rm so}_2\oplus {\rm so}_2)$.
\end{abstract}

\vfill

\newpage

\centerline{\sc 1. Introduction}
\medskip

\noindent
   In the classical case, the embedding $SO(n)\subset U(n)$ is of great
importance in the theory of Riemannian spaces and for
group theoretical approach to some physical problems. In the frame of
Drinfeld--Jimbo quantum groups we cannot construct the corresponding
embeddings. In this frame we also cannot construct the quantum
algebras $U_q({\rm {so}}_{n,1})$ and introduce Gel'fand--Tsetlin bases in
spaces of irreducible representations of $U_q({\rm {so}}_n)$. To remove
these defects the
new $q$-deformation of the universal enveloping algebra
$U({\rm {so}}(n,{\bf C}))$
was defined in [1] (see also [2]). We denote it by $U'_q({\rm {so}}_n)$.
This $q$-deformed algebra
allows the embedding $U'_q({\rm {so}}_{n-1})
\subset U'_q({\rm {so}}_n)$ and, therefore, we can introduce
Gel'fand--Tsetlin bases. It was shown by Noumi [3]
that this algebra can be embedded into $U_q({\rm {sl}}_n)$.
The last fact make the algebra $U'_q({\rm {so}}_n)$ very attractive since
the pair $U'_q({\rm {so}}_n)\subset U_q({\rm sl}_n)$ is of great
importance in mathematics and theoretical physics. Of course, for
applications we must have good developed representation theory
 of $U'_q({\rm {so}}_n)$. Representations of this algebra
for $q$ not a root of unity where studied in several papers (see [1, 2]).
For $q$ a root of unity only representations of the algebra
$U'_q({\rm {so}}_3)$ were investigated [4, 5].

The algebra $U'_q({\rm {so}}_n)$ has real forms $U'_q({\rm {so}}_{m,r})$.
The representation theory of the algebras $U'_q({\rm {so}}_{2,1})$ and
$U'_q({\rm {so}}_{3,1})$ is developed in [2].
The aim of this paper is to construct irreducible representations of
$U'_q({\rm {so}}_{2,2})$ when $q$ is a positive number.
In order to construct representations of the algebra $U'_q({\rm {so}}_{2,2})$
we derive formulas for irreducible finite dimensional representations of
$U'_q({\rm {so}}_{4})$ (see [1] and [6])
in the basis corresponding to the reduction from
$U'_q({\rm {so}}_{4})$ to $U({\rm so}_2\oplus {\rm so}_2)$. Then we
"analytically continue" these formulas to obtain infinite dimensional
representations of $U'_q({\rm {so}}_{2,2})$. In this way we obtain the
series of representations $T^\epsilon _{bc}$ given by two complex
parameters. Theorem 1 below describe when these representations are
irreducible. Then we study reducible representations $T^\epsilon _{bc}$ and
separate irreducible constituents from them. In this way we obtain Theorem 2
describing all irreducible representations which can be obtained from the
representations $T^\epsilon _{bc}$. Then in Theorem 3
we separate all $*$-representations
from the set of irreducible representations of Theorem 2.
Comparing Theorems 2 and 3 with the representation theory of the Lie algebra
${\rm so}_{2,2}$, we see that there exists big difference between
representation theories of $U'_q({\rm {so}}_{2,2})$ and so$_{2,2}$. The
same can be said about representation theories of $U'_q({\rm {so}}_{2,2})$
and $U'_q({\rm {so}}_{3,1})$. These differences are emphasized in Conclusions
at the end of the paper.
\bigskip

\centerline{{\sc 2. The algebra} $U'_q({\rm {so}}_{2,2})$}
\medskip

\noindent
Drinfeld [7] and Jimbo [8] defined $q$-deformed (quantum)
algebras $U_q(g)$ for all simple complex Lie algebras $g$ by means of
Cartan subalgebras and root subspaces (see also [9]).
However, these approaches do not give a satisfactory presentation of the
quantum algebra $U_q({\rm so}_n)$
from point of view of some problems of quantum physics and representation
theory. In fact, they admit the inclusion
$U_q({\rm so}_n)\supset U_q({\rm so}_{n-2})$
and do not admit the inclusion
$$U_q({\rm so}_n)\supset U_q({\rm so}_{n-1}). \eqno(1)$$
This is why we cannot construct the quantum algebra $U_q({\rm so}_{n,1})$
in the frame of these approaches and cannot construct Gel'fand--Tsetlin
bases in the representation spaces. In order to obtain inclusion (1) it was
proposed in [1] another $q$-deformation of the
classical universal enveloping algebra $U({\rm so} (n,{\bf C}))$. The
classical
algebra $U({\rm so} (n,{\bf C}))$ is generated by the elements $I_{i,i-1}$,
$i=2,3,\ \cdots ,n$, that satisfy the relations
$$I_{i,i-1}I^2_{i+1,i}-2I_{i+1,i}I_{i,i-1}I_{i+1,i}+I^2_{i+1,i}I_{i,i-1}
=-I_{i,i-1}, \eqno(2)$$
$$I^2_{i,i-1}I_{i+1,i}-2I_{i,i-1}I_{i+1,i}I_{i,i-1}+I_{i+1,i}I^2_{i,i-1}
=-I_{i+1,i}, \eqno(3)$$
$$[I_{i,i-1},I_{j,j-1}]=0,\ \ \vert i-j\vert >1. \eqno(4)$$
They follow from the well-known commutation relations for the generators
$I_{ij}$ of the Lie algebra ${\rm so} (n,{\bf C})$ (see the paper by
Gel'fand and Tsetlin [10]).

   The approach of the paper [1] to the $q$-deformed orthogonal
algebra consists in
a $q$-deformation of the associative algebra $U({\rm so} (n,{\bf C}))$ by
deforming relations (2)-(4). The $q$-deformed relations are of the form
$$I_{i,i-1}I_{i+1,i}^2-aI_{i+1,i}I_{i,i-1}I_{i+1,i}
+I^2_{i+1,i}I_{i,i-1}=-I_{i,i-1}, \eqno(5)$$
$$I^2_{i,i-1}I_{i+1,i}-aI_{i,i-1}I_{i+1,i}I_{i,i-1}
+I_{i+1,i}I^2_{i,i-1}=-I_{i+1,i}, \eqno(6)$$
$$[I_{i,i-1},I_{j,j-1}]=0,\ \ \vert i-j\vert >1, \eqno(7)$$
where $a=q+q^{-1}=(q^2-q^{-2})/(q-q^{-1})$
and $[\cdot ,\cdot ]$ denotes the usual commutator. Obviously, in the
limit $q\to 1$ formulas (5)--(7) give relations (2)--(4). Remark that
relations (5) and (6) differ from the $q$-deformed Serre relations in the
approach of Jimbo and Drinfeld to quantum orthogonal algebras by
appearance of nonzero right hand side and by a possibility of reduction (1).
Below, by the algebra $U'_q({\rm so}_n)$ we mean the $q$-deformed
algebra defined by formulas (5)--(7). Unfortunately, the algebra
$U'_q({\rm so}_n)$ does not have a Hopf algebra structure. But it can be
embedded into the Hopf algebra $U_q({\rm sl}_n)$ and is a Hopf ideal
in $U_q({\rm sl}_n)$ [3].

   As in the classical case, the $q$-algebras $U'_q({\rm so}_3)$ and
$U'_q({\rm so}_4)$ can also be described in terms of bilinear
relations
($q$-commutators). In fact, defining the algebra $U'_q({\rm so}_3)$
by relations
(5)--(7) we have only two generators $I_{21}$ and $I_{32}$. However, we
can define the third element $I_{31}$ according to the formula
$$I_{31}=q^{1/2}I_{21}I_{32}-q^{-1/2}I_{32}I_{21} \eqno (8) $$
(see [2]).
Then by the algebra $U'_q({\rm so}_3)$ we mean the associative algebra
generated by the elements $I_{21},$ $I_{32}$ and $I_{31}$ which satisfy
the relations
$$q^{1/2}I_{21}I_{32}-q^{-1/2}I_{32}I_{21}=I_{31},\eqno (9) $$
$$q^{1/2}I_{31}I_{21}-q^{-1/2}I_{21}I_{31}=I_{32},\eqno (10) $$
$$q^{1/2}I_{32}I_{31}-q^{-1/2}I_{31}I_{32}=I_{21}.\eqno (11) $$

   It is clear that if the generators $I_{21},$  $I_{32}$ and $I_{31}$ satisfy
the relations (9)--(11), then the pair $I_{21}$ and $I_{32}$ satisfies the
trilinear relations (5) and (6). Remark that the algebra given by
relations (9)--(11) coincides with the cyclically symmetric
algebra defined in [11].

   The $q$-deformed algebra $U'_q({\rm so}_4)$
is generated by the generators $I_{21}$, $I_{32}$ and $I_{43}$.
Moreover, for the first two generators everything, said above
$U'_q({\rm so}_3)$, is true. Thus, the inclusion
$U'_q({\rm so}_3) \subset U'_q({\rm so}_4)$
takes place. The generators $I_{21}$ and $I_{43}$ mutually commute
(see relation (7))
and the pair $I_{32},$ $I_{43}$ in turn must satisfy relations (5) and
(6). Again, $U'_q({\rm so}_4)$ can be also given in terms of bilinear
$q$-commutators. Namely, we can add to the triple of generators $I_{21},$
$I_{32}$ and $I_{43}$ the element $I_{31}$ from (8) and the elements
$I_{42}$, $I_{41}$ defined as
$$I_{42}=q^{1/2}I_{32}I_{43}-q^{-1/2}I_{43}I_{32}, $$
$$I_{41}=q^{1/2}I_{31}I_{43}-q^{-1/2}I_{43}I_{31}=q^{1/2}I_{21}I_{42}-q^{-1/2}
I_{42}I_{21}.  $$

    Various real forms of the algebra $U'_q({\rm so}_4)$ are obtained
by introducing corresponding $*$-structures (antilinear antihomomorphisms).
The real form $U'_q({\rm so}_{2,2})$ is defined by the $*$-structure
$$
I^*_{21}=-I_{21},\ \ \ I^*_{32}=I_{32},\ \ \ I^*_{34}=-I_{34}.
$$
The relations
$I^*_{21}=-I_{21}$, $I^*_{32}=-I_{32}$, $I^*_{34}=I_{34}$
determine the real form $U'_q({\rm so}_{3,1})$.

Everywhere below we assume that $q$ is a positive real number.
\bigskip

\centerline{\sc 3. Representations of the algebra $U'_q({\rm so}_4)$}
\medskip

\noindent
Let us describe irreducible finite dimensional representations of
$U'_q({\rm so}_4)$ when $q$ is not a root of unity.
They are described by using representations of $U'_q({\rm so}_3)$.

Irreducible finite dimensional representations of $U'_q({\rm so}_3)$
are given by integral or half-integral nonnegative
number $l$. We denote these representations by $T_l$. The carrier space of
the representation $T_l$ has the orthonormal basis $\{ |m\rangle$,
$m=l,l-1,\cdots ,-l\}$, and the operators $T_l(I_{21})$ and $T_l(I_{32})$
act upon this basis as
$$T_l(I_{21})|m\rangle ={\rm i}[m]|m\rangle ,\eqno (12)$$
$$T_l(I_{32})|m\rangle =d(m)([l-m][l+m+1])^{1/2}|m+1\rangle
-d(m-1)([l-m+1][l+m])^{1/2}|m-1\rangle ,\eqno (13)$$
where
$$d(m)=([m][m+1]/[2m][2m+2])^{1/2}$$
and $[a]$ denotes a $q$-number defined by
$$[a]=(q^{a}-q^{-a})/(q-q^{-1}).$$
If $q$ is positive then these representations exhaust all irreducible
finite dimensional representations of $U'_q({\rm so}_3)$. For other values
of $q$ there exist irreducible finite dimensional representations
which are not equivalent to these representations (see, for example, [4])
but we shall not need them.

As in the case of the Lie group $SO(4)$, finite dimensional
irreducible representations $T_{rs}$ of the $q$-deformed algebra
$U'_q({\rm so}_4)$ are given by two integral or half-integral numbers
$r$ and $s$
such that $r\ge |s|\ge 0$ (see [6]). Restriction of $T_{rs}$
onto the subalgebra $U'_q({\rm so}_3)$ decomposes into the sum of the
irreducible representations $T_l$ of this subalgebra for which $l=|s|,
|s|+1,\cdots ,r$. Uniting the bases of the subspaces
of the irreducible representations
$T_l$ of $U'_q({\rm so}_3)$ we obtain the basis of the carrier space
$V_{rs}$ of the representation $T_{rs}$ of $U'_q({\rm so}_4)$. Thus, the
corresponding orthonormal basis of $V_{rs}$ consists of the vectors
$$|l,m\rangle ,\ \ \ |s|\le l\le r,\ \ \ m=-l,-l+1,\cdots ,l.$$
The operator $T_{rs}(I_{43})$ acts upon these vectors by the formula
$$T_{rs}(I_{43})|l,m\rangle ={\rm i}{[r+1][s][m]\over [l][l+1]}|l,m\rangle
\qquad\qquad\qquad\qquad\qquad\qquad\qquad$$
$$+\left( {[r-l][l+s+1][l-s+1][l+m+1][l-m+1]\over [r+l+2]^{-1}[l+1]^2[2l+1]
[2l+3]}\right) ^{1/2}|l+1,m\rangle $$
$$-\left( {[r+l+1][l+s][l-s][l+m][l-m]\over [r-l+1][l]^2[2l-1][2l+1]}\right)
^{1/2}|l-1,m\rangle ,\eqno (14)$$
where numbers in the square brackets are $q$-numbers.
The operators $T_{rs}(I_{21})$ and $T_{rs}(I_{32})$ act upon the basis
vectors by formulas (12) and (13). Formulas (12), (13) and (14)
completely determine the representation $T_{rs}$.
\bigskip

\centerline{{\sc 4. Diagonalization of the operator} $T_{rs}(I_{43})$}
\medskip

\noindent
We shall need the representations $T_{rs}$ in another form. To obtain it
we diagonalize the operator $T_{rs}(I_{43})$. In the next section
this result is used for obtaining representations $T_{rs}$ in
the bases corresponding to restriction upon the subalgebra
$U'_q({\rm so}_2)+U'_q({\rm so}_2)$. It is more convenient to deal with the
selfadjoint operator $L=-{\rm i}T_{rs}(I_{43})$, ${\rm i}=\sqrt {-1}$.
Replacing the
vectors $\vert l,m\rangle $ by $\vert l,m\rangle '={\rm i}^{-l}\vert l,m
\rangle $ we obtain that $L$ acts upon the vectors $\vert l,m\rangle '$
by formula (14) in which the sign -- of the third summand is replaced
by + and the first summand is multipled by $-{\rm i}$.

   The space $V_{rs}$ can be decomposed into the sum
$V_{rs}=\sum ^r_{m=-s} \oplus \ V_m ,$
where $V_m$ is spanned by the vectors $\vert l,m\rangle $ with fixed $m$.
Let us find the spectrum and the eigenvectors
$$\vert x,m \rangle '=\sum ^r_{l=k} P_{l-k}(x)\vert l,m\rangle ,\ \ \
k=\max \ (\vert m\vert ,\ \vert s\vert ) \eqno(15)$$
of the operator $L$ on the subspace $V_m$:
$$L\vert x,m \rangle '=[x]\vert x,m \rangle ',\eqno (16)$$
where $[x]$ is a $q$-number. Formula (14) is
symmetric with respect to permutation of $s$ and $m$ and to change
of signs at $m$ and $s$. Therefore, we may assume, without loss of
generality, that $s$ and $m$ are positive and that $s\ge m$.

   Substituting expression (15) for $\vert x,m \rangle '$ into (16) and
acting by $L$ upon $|l,m\rangle$ we easily find that vector
(15) is an eigenvector of $L$ with the eigenvalue $[x]$
if $P_{l-k}$ satisfy the recurrence relation
$$\Biggl( {[u][n+2s+1][n+1][n+s+m+1][n+s-m+1]\over
[r+n+s+2]^{-1}[n+s+1]^2[2n+2s+1][2n+2s+3]}\Biggr) ^{1/2} P_{n+1}(x)$$
$$+\Biggl( {[r+n+s+1][r-n-s+1][n+2s][n][n+s+m]\over [n+s-m]^{-1}
[n+s]^2[2n+2s-1][2n+2s+1]}\Biggr) ^{1/2}P_{n-1}(x)$$
$$+{[r+1][s][m]\over [n+s][n+s+1]}P_n(x)=[x]P_n(x) \eqno(17)$$
(here $u=r-n-s$, $n=l-k$) and the initial conditions $P_0(x)=1$, $P_{-1}(x)=0$.

   Making in (17) the substitution
$$P_n(x)=-q^c\Biggl( {[n+2s]![n+s+m]![2n+2s+1]\over
[n]![n+s-m]![r-n-s]![r+n+s+1]!}\Biggr) ^{1/2}P_n'(x)$$
where $c=s+m-r$, we reduce (17) to the recurrence relation
$${(1-Q^{n+2s+1})(1-Q^{n+s+m+1})(1-Q^{n-r+s})(1+Q^{n+s+1})\over
(1-Q^{2n+2s+1})(1-Q^{2n+2s+2})}P_{n+1}'(x)$$
$$-{Q^{s+m-r}(1-Q^n)(1-Q^{r+n+s+1})(1+Q^{n+s})(1-Q^{n+s-m})\over
(1-Q^{2n+2s+1})(1-Q^{2n+2s})}P_{n-1}'(x)$$
$$-{Q^{n-r+s}(1-Q^{r+1})(1-Q^m)(1-Q^s)\over (1-Q^{n+s})(1-Q^{n+s+1})}
P_n'(x)={q-q^{-1}\over Q^{(r-s-m)/2}}[x]P_n'(x),$$
where $Q:=q^2$.
Comparing this formula with recurrence relation (7.5.2) from the book of
Gasper and Rahman [13] for $q$-Racah polynomials
$$R_n(\mu (y);\ \alpha ,\beta , \gamma ,\delta \vert Q)
={}_4\varphi _3 \left(\matrix{Q^{-y},Q^{y+1}\gamma \delta ,Q^{-n},
Q^{n+1}\alpha \beta \cr
\alpha Q,\ \ \beta \delta Q,\ \ \gamma Q}
;\ Q,Q\right) $$
(here ${}_4\varphi _3$ is a basic hypergeometric function which can be
found in [13]) at
$$\alpha =\beta =-Q^s,\ \gamma =Q^{s+m},\ \delta =-Q^{-r-1}, \eqno(18)$$
after cumbersome transformations we conclude that
$$P'_n(x)=R_n(\mu (y);\ \alpha ,\beta , \gamma ,\delta \vert Q),$$
where $\alpha ,\ \beta ,\ \gamma ,\
\delta $ are given by formulas (18) and
$x=(r-s-m)-2y$.
Thus, the polynomials $P_n(x)$ from (17) normalized by the condition
$P_0(x)=1$ are of the form
$$P_n(x)=N^{1/2}R_n(\mu (y);\ -Q^s,-Q^s,Q^{s+m},-Q^{-r-1}\vert Q), \eqno(19)$$
$$N={[n+2s]![n+s+m]![2n+2s+1][s-m]![r-s]![r+s+1]!\over
[n]![n+s-m]![r-n-s]![r+n+s+1]![2s]![s+m]![2s+1]},$$
where $x=(r-s-m)-2y$. The variable $y$ takes the values
$0,1,2,\ \cdots ,\ r-s$. Therefore, the spectrum of $L$ on
the subspace $V_m$ consists of the points
$$[r-s-m],\ [r-s-2-m],\ [r-s-4-m],\ \cdots ,\ [-(r-s)-m]. \eqno(20)$$
The corresponding eigenvectors are determined by formulas (15) and
(19). The orthogonality relation for the polynomials $P_n(x)$ follows
from the orthogonality of $q$-Racah polynomials (see [13]) and
is of the form
$$\sum ^{r-s}_{y=0} P_n(x)P_k(x)W(x)=\delta _{nk}.\eqno (21) $$
Here $W(x)$ is equal to the expression
$${[4y+2k-2r][2y+2k-2r-2]!![2y+2s]!![2r-2y]!![r-m-y]!
\over [2y+2k-2r][y+k-r-1]![y+s]!
[2y+2m]!![r-y]![r-s-y]![y]!} $$
$$\times [y+m]![k+y]![2s+1]!!([s]!)^2[r-s]([2s]!![s-m]![k]![r+s+1]!)^{-1},$$
where $k=s+m$, $[n]!=[n][n-1]\cdots [1]$ and $[n]!!=[n][n-2][n-4]\cdots
[1]$ or $[2]$.

Formula (21) shows that vectors (15) are not normalized. The vectors
$|x,m\rangle =W(x)^{1/2}|x,m\rangle '$
are normal and due to formula (16) we have
$$T_{rs}(I_{43})|x,m\rangle ={\rm i}[x]|x,m\rangle .\eqno (22)$$

    Joining spectra (20) for all subspaces $V_m$, we obtain the spectrum of
the operator $T_{rs}$, and therefore the spectrum of the operator
$T_{rs}(I_{43})$.
\bigskip

\centerline{{\sc 5. Representations} $T_{rs}$ {\sc in the basis} $|x,m\rangle$}
\medskip

\noindent
The operator $T_{rs}(I_{43})$ acts upon the basis vectors
$|x,m\rangle$ by formula (22). It is clear from formulas (13) and (15)
that
$$T_{rs}(I_{21})|x,m\rangle ={\rm i}[m]|x,m\rangle .\eqno (23)$$
Thus, to have the representation $T_{rs}$ in the basis $|x,m\rangle$,
we must find the action formula for the operator $T_{rs}(I_{32})$
upon this basis.

Since
$$|x,m\rangle =\sum _{l=s}^r P^m_{l-s}(x)|l,m\rangle ,\eqno (24)$$
with $P^m_{l-s}(x)=W(x)^{1/2}P_{l-s}(x)$, then due to formula (13) we have
$$T_{rs}(I_{32})|x,m\rangle =d(m)\sum ^r_{l=s}P^m_{l-s}(x)([l-m][l+m+1])^{1/2}
|l,m+1\rangle$$
$$-d(m-1)\sum ^r_{l=s}P^m_{l-s}(x)([l-m+1][l+m])^{1/2}|l,m-1\rangle .
\eqno (25)$$
Applying to $([l-m][l+m+1])^{1/2}P^m_{l-s}(x)$ recurrence relation
(7.2.14) of [13] with
$$a=Q^{m-r-1},\ \ \ b=-Q^s,\ \ \ c=d=-Q^m,\ \ \ n=(r-s-m-x)/2,\ \ \
j=l-m$$
and using the equalities
$$[2x]/[x]=q^{x}+q^{-x},\ \ \ \ (q^{a+b}\pm q^{-(a+b)})
(q^{a-b}\mp q^{-(a-b)})=[2a]\mp [2b],$$
after some calculations we obtain for the first summand of the right hand
side of (25) the expression
$$d(m)d(x-1)\{ ([r+1]+[s-m+x-1])([r+1]+[s+m-x+1])\} ^{1/2}
\sum ^r_{l=s}P^{m+1}_{l-s}(x-1)|l,m+1\rangle$$
$$-d(m)d(x)\{ ([r+1]-[s+m+x+1])([r+1]-[s-m-x-1])\} ^{1/2}
\sum ^r_{l=s}P^{m+1}_{l-s}(x+1)|l,m+1\rangle . \eqno (26)$$

   To transform the second summand on the right hand side of (25) we
apply to the basic hypergeometric function $_4\varphi _3$ from the
expression for $P^m_{l-s}(x)$ the transformation
$${}_4\varphi _3 \!\! \left(\matrix{Q^{-N}, \alpha ,  \beta , \gamma \cr
\delta ,\ \ \sigma ,\ \ \rho }
; Q,Q\right)
={(\sigma /\alpha ;Q)_N(\rho /\alpha ;Q)_N\over (\sigma ;Q)_N(\rho ;Q)_N
\alpha ^{-N}}\
_4\varphi _3 \!\! \left(\matrix{Q^{-N},\alpha ,\  \delta /\beta ,\
\delta /\gamma \cr
\delta , \alpha Q^{1-N}/\sigma , \alpha Q^{1-N}/\rho }
; Q,Q\right) $$
(see [13]), where $N=l-s$ and
$$\alpha =Q^{l+s+1},\ \beta =-Q^{(s-r+m-x)/2},\ \gamma =Q^{(s-r+m+x)/2},$$
$$
\delta =Q^{s-r},\ \sigma =-Q^{s+1},\ \rho =Q^{s+m+1}.$$
Here $(a;Q)_n=(1-Q)(1-aQ)(1-aQ^2)\cdots (1-aQ^{n-1})$. Now we apply to
$([l-m+1][l+m])^{1/2}P^m_{l-s}(x)$ the same recurence relation (7.2.14)
of [13] with
$$a=Q^{-(r+m+1)},\ \ \ b=-Q^{-s},\ \ \ c=d=-Q^m,\ \ \ n=(r-s+m+x)/2,\ \ \
j=l+m.$$
Then the second summand of the right hand side of (25) takes the form
$$d(m-1)d(x)\{([r+1]+[s+m-x-1])([r+1]+[s-m+x+1])\} ^{\frac 12}
\sum ^r_{l=s} \! P^{m-1}_{l-s}(x+1)|l,m-1\rangle$$
$$-d(m-1)d(x-1)\{([r+1]-[s-m-x+1])([r+1]-[s+m+x-1])\} ^{\frac 12}
\sum ^r_{l=s} \! P^{m-1}_{l-s}(x-1)|l,m-1\rangle  \eqno (27)$$

   We substitute expressions (26) and (27) into (25) and take into account
formula (24). As a result, we find that the operator $T_{rs}(I_{32})$ acts
upon the vectors $|x,m\rangle$ as
$$T_{rs}(I_{32})|x,m\rangle = $$
$$
=d(m)d(x-1)\{ ([r+1]+[s-m+x-1])([r+1]+[s+m-x+1])
\} ^{1/2}|x-1,m+1\rangle$$
$$-d(m)d(x)\{ ([r+1]-[s+m+x+1])([r+1]-[s-m-x-1])
\} ^{1/2}|x+1,m+1\rangle$$
$$+d(m-1)d(x-1)\{ ([r+1]-[s-m-x+1])([r+1]-[s+m+x-1])
\} ^{1/2}|x-1,m-1\rangle$$
$$-d(m-1)d(x)\{ ([r+1]+[s+m-x-1])([r+1]+[s-m+x+1])
\} ^{1/2}|x+1,m-1\rangle . \eqno (28)$$
Now we completely determined representations $T_{rs}$ of $U'_q({\rm {so}}_4)$
with respect to the basis corresponding to reduction onto the subalgebra
$U'_q({\rm {so}}_2)+U'_q({\rm {so}}_2)$.
\bigskip

\centerline{{\sc 6. Infinite dimensional representations of}
$U'_q({\rm {so}}_{2,2})$}
\medskip

\noindent
Let us first define infinite dimensional linear representations of
$U'_q({\rm so}_{2,2})$. By a linear representation $T$ of the algebra
$U'_q({\rm so}_{2,2})$ we mean a homomorphism of $U'_q({\rm so}_{2,2})$
into the algebra of linear operators (bounded or unbounded) on a
Hilbert space $H$, defined on an everywhere dense invariant subspace $D$,
such that
\medskip

(a) the operators $T(I_{21})$ and $T(I_{43})$ can be simulteniously
diagonalized,

(b) eigenvalues of $T(I_{21})$ and $T(I_{43})$ have finite multiplicities,

(c) eigenvectors of $T(I_{21})$ and $T(I_{43})$ belong to $D$.
\medskip

A representation $T$ of $U'_q({\rm so}_{2,2})$ is called a $*$-representation
if the operators $T(I_{21})$, $T(I_{32})$ and $T(I_{43})$ satisfy on $D$
the relations
$$
T(I_{21})^*=-T(I_{21}),\ \ \
T(I_{32})^*=T(I_{32}),\ \ \
T(I_{43})^*=-T(I_{43}).
$$

As in the case of representations of compact and noncompact real
Lie groups, by making use of analytical continuation in parameters giving
representations we can obtain infinite dimensional representations of the
$q$-deformed algebra $U'_q({\rm {so}}_{2,2})$ from the representations
$T_{rs}$ of $U'_q({\rm {so}}_4)$. In this way, we obtain the representations
$T^\epsilon _{\sigma \tau }$, $\sigma \in {\bf C}$, $\tau \in {\bf C}$,
$\epsilon \in \{ 0,\ 1\}$, of $U'_q({\rm {so}}_{2,2})$ which act on the
Hilbert spaces $H_\epsilon$ with the orthonormal basis
$$|x,m\rangle ,\ \ \ x\in {\frac 12}{\bf Z},\ \ \
m\in {\frac 12}{\bf Z},\ \ \
x+m\equiv \epsilon \ ({\rm mod}\ 2).$$
The operators $T^\epsilon _{\sigma \tau }(I_{21})$ and
$T^\epsilon _{\sigma \tau }(I_{43})$ act upon these basis vectors by
formulas (22) and (23). For the operator $T^\epsilon _{\sigma \tau }(I_{32})$
we have
$$T^\epsilon _{\sigma \tau }(I_{32}) |x,m\rangle$$
$$=d(m)d(x-1)\{ ([\sigma +1]+[\tau -m+x-1])([\sigma +1]+[\tau +m-x+1])
\} ^{\frac 12}|x-1,m+1\rangle$$
$$-d(m)d(x)\{ ([\sigma +1]-[\tau +m+x+1])([\sigma +1]-[\tau -m-x-1])
\} ^{\frac 12}|x+1,m+1\rangle$$
$$+d(m-1)d(x-1)\{ ([\sigma +1]-[\tau -m-x+1])([\sigma +1]-[\tau +m+x-1])
\} ^{\frac 12}|x-1,m-1\rangle$$
$$-d(m-1)d(x)\{ ([\sigma +1]+[\tau +m-x-1])([\sigma +1]+[\tau -m+x+1])
\} ^{\frac 12}|x+1,m-1\rangle . $$

It is more convenient to write down the last formula in the form
$$T^\epsilon _{\sigma \tau }(I_{32}) |x,m\rangle$$
$$=\left( {[\sigma -\tau +m-x+2][\sigma -\tau -m+x]
[(\sigma +\tau +m-x+2)/2]\over
[(\sigma -\tau +m-x+2)/2][(\sigma -\tau -m+x)/2]
[(\sigma +\tau -m+x)/2]^{-1}}\right) ^{1/2}$$
$$\times d(m)d(x-1)|x-1,m+1\rangle$$
$$-\left( {[\sigma +\tau +m+x+2][\sigma +\tau -m-x]
[(\sigma -\tau +m+x+2)/2]\over
[(\sigma +\tau +m+x+2)/2][(\sigma +\tau -m-x)/2]
[(\sigma -\tau -m-x)/2]^{-1}}\right) ^{1/2}$$
$$\times d(m)d(x)|x+1,m+1\rangle$$
$$+\left( {[\sigma +\tau -m-x+2][\sigma +\tau +m+x]
[(\sigma -\tau -m-x+2)/2]\over
[(\sigma +\tau -m-x+2)/2][(\sigma +\tau +m+x)/2]
[(\sigma -\tau +m+x)/2]^{-1}}\right) ^{1/2}$$
$$\times d(m-1)d(x-1)|x-1,m-1\rangle$$
$$-\left( {[\sigma -\tau -m+x+2][\sigma -\tau +m-x]
[(\sigma +\tau -m+x+2)/2]\over
[(\sigma -\tau -m+x+2)/2][(\sigma -\tau +m-x)/2]
[(\sigma +\tau +m-x)/2]^{-1}}\right) ^{1/2}$$
$$\times d(m-1)d(x)|x+1,m-1\rangle .\eqno (29)$$
In every summand here there are two expressions of the form
$[\sigma -\tau -m+x]/[(\sigma -\tau -m+x)/2]$.
This expression is equal to
$q^{(\sigma -\tau -m+x)/4}+q^{-(\sigma -\tau -m+x)/4}.$

Introducing the notation
$$
\gamma =\sigma +\tau +2,\ \ \ \delta =\sigma -\tau +2,\ \ \
M=x+m,\ \ \ N=m-x,
$$
where $M$ and $N$ are even if $\epsilon =0$ and odd if $\epsilon =1$,
and denoting the basis elements $|x,m\rangle $ by $|M,N\rangle '$,
we obtain after multiplication of $|M,N\rangle '$ by the
appropriate factors that
$$
T^\epsilon _{\sigma \tau }(I_{32}) |M,N\rangle =
([(M+N)/2][(M-N)/2])/([M+N][M-N]) \qquad\qquad\qquad\qquad\qquad$$
$$\qquad\qquad
\times \{ (q^{(\delta +N)/2}+q^{-(\delta +N)/2})[(\gamma +N)/2] |M,N+2\rangle $$
$$\qquad\qquad
-(q^{(\gamma +M)/2}+q^{-(\gamma +M)/2})[(\delta +M)/2] |M+2,N\rangle $$
$$\qquad\qquad
+(q^{(\gamma -M)/2}+q^{-(\gamma -M)/2})[(\delta -M)/2] |M-2,N\rangle $$
$$\qquad\qquad
-(q^{(\delta -N)/2}+q^{-(\delta -N)/2})[(\gamma -N)/2] |M,N-2\rangle \} ,
$$
where $|M,N\rangle $ are the basis elements $|M,N\rangle '$ with the
appropriate factors.

Setting
$$
k=M/2,\ \ \ l=N/2,\ \ \ c=\gamma /2,\ \ \ b=\delta /2
$$
and denoting the basis elements $|M,N\rangle $ by $|k,l\rangle $ and the
representations $T^\epsilon _{\sigma \tau }$ by $T^\epsilon _{bc}$
we obtain the representations in the form
$$
T^\epsilon _{bc}(I_{32})|k,l\rangle =([k+l][k-l])/([2(k+l)][2(k-l)])
\qquad\qquad\qquad\qquad\qquad$$
$$
\times \{ (q^{l+b}+q^{-l-b})[l+c] |k,l+1\rangle -
(q^{k+c}+q^{-k-c})[k+b] |k+1,l\rangle $$
$$
+(q^{-k+c}+q^{k-c})[b-k] |k-1,l\rangle
-(q^{b-l}+q^{l-b})[c-l] |k,l-1\rangle \} , \eqno (30) $$
$$
T^\epsilon _{bc}(I_{21})|k,l\rangle =\sqrt {-1} [k+l]|k,l\rangle ,\ \ \ \ \
T^\epsilon _{bc}(I_{43})|k,l\rangle =\sqrt {-1} [k-l]|k,l\rangle .\eqno (31)
$$
Note that the basis consists of the vectors $|k,l\rangle $,
where $k$ and $l$ are integral if $\epsilon =0$ and half-integral (half of
an odd integer) if $\epsilon =1$. We consider that the invariant everywhere
dense subspace $D_\epsilon$ of $H_\epsilon$ coincides with the span of all
vectors $| k,l\rangle$ from $H_\epsilon$.

Remark that the operators $T^\epsilon _{bc}(I_{21})$ and
$T^\epsilon _{bc}(I_{43})$ are unbounded and the operator
$T^\epsilon _{bc}(I_{32})$ is bounded.

Our aim is to study representations $T^\epsilon _{bc}$ of
$U'_q({\rm {so}}_{2,2})$. We say that $T^\epsilon _{bc}$ is irreducible if
it is algebraically irreducible on the subspace $D_\epsilon$.
\medskip

\noindent
{\bf Theorem 1.} {\it The representation $T^\epsilon _{bc}$ is irreducible
if and only if no of the numbers $b$ and $c$ coincides with any of the
numbers $n$, $n+{\rm i}\pi r/2h$, $n,r\in {\bf Z}$, where $h$ is defined by
$q=\exp h$.}
\medskip

\noindent {\it Proof} is given as in the case of representations of
semisimple Lie algebras (see, for example, [14], Chapter 7).
\medskip

There exist equivalence relations in the set of representations
$T^\epsilon _{bc}$. One type of equivalences appears because of the
periodicity of the function $w(z)=[z]$, where $[z]=(q^z-q^{-z})/(q-q^{-1})$.
If $q=\exp h$, then the function $w(z)$ is periodic with period
$2\pi {\rm i}/h$. Therefore, it follows from (30) and (31) that
$$
T^\epsilon _{bc}=T^\epsilon _{b+2\pi {\rm i}/h,c} =
T^\epsilon _{b,c+2\pi {\rm i}/h}. \eqno (32)
$$
For the function $w(z)$ we also have $w(z)=-w(z+\pi {\rm i}/h)$.
For this reason, replacement of $b$ (respectively, of $c$) by
$b+\pi {\rm i}/h$ (respectively, by $c+\pi {\rm i}/h$)
in the relation (30) changes only signs near vectors,
and the representations $T^\epsilon _{bc}$ and
$T^\epsilon _{b+\pi {\rm i}/h,c}$ (respectively,
$T^\epsilon _{b,c+\pi {\rm i}/h}$) are equivalent and the equivalence
operator is diagonal with respect to the basis $\{ | k,l\rangle \}$
with numbers $\pm 1$ on the main diagonal. Thus,
$$
T^\epsilon _{bc}\sim T^\epsilon _{b+\pi {\rm i}/h,c}\sim
T^\epsilon _{b,c+\pi {\rm i}/h}. \eqno (33)
$$
If the representation $T^\epsilon _{bc}$ is irreducible, then we also have
the equivalences
$$
T^\epsilon _{bc} \sim T^\epsilon _{b,-c+1} \sim T^\epsilon _{-b+1,c}.
\eqno (34)
$$
The equivalence operators are diagonal in the basis $\{ | k,l\rangle \}$
and their matrix elements are calculated in the same way as in the
case of representations of Lie algebras (see, for example, [15],
Sect. 6.4.4).

Using the method of [15], Sect. 6.4, it is easy to prove that any
equivalence relation between irreducible representations in the set of
representations $T^\epsilon _{bc}$ is a composition of equivalence relations
given above.

Taking into account the equivalences (32)--(34),
everywhere below we assume (without loosing the generality) that
$0\le {\rm Im}\ b <\pi{\rm i}/h$, $0\le {\rm Im}\ c <\pi{\rm i}/h$,
Re$\ b \ge 1/2$ and Re$\ c \ge 1/2$.

\bigskip

\centerline{{\sc 7. Irreducible subrepresentations of} $T^\epsilon _{bc}$}
\medskip

In order to find irreducible constituents of reducible representations
$T^\epsilon _{bc}$ we reduce these representations to the form (29). If
the representation $T^\epsilon _{bc}$ is irreducible, then for the
operator $T^\epsilon _{bc}(I_{32})$ we have
$$
T^\epsilon _{bc}(I_{32})| k,l\rangle '=([k+l][k-l])/([2(k+l)][2(k-l)])$$
$$
\times \{ ((q^{l+b}+q^{-l-b})(q^{b-l-1}+q^{l-b+1})[l+c][l-c+1])^{1/2}
| k,l+1\rangle ' $$
$$
- ((q^{k+c}+q^{-k-c})(q^{c-k-1}+q^{k-c+1})[k+b][k-b+1])^{1/2}
| k+1,l\rangle ' $$
$$
- ((q^{k+c-1}+q^{-k-c+1})(q^{c-k}+q^{k-c})[k-b][k+b-1])^{1/2}
| k-1,l\rangle ' $$
$$
+((q^{l+b-1}+q^{-l-b+1})(q^{b-l}+q^{l-b})[l-c][l+c-1])^{1/2}
| k,l-1\rangle '\} , \eqno (35)
$$
where the vectors $| k,l\rangle '$ are obtained from $| k,l\rangle $ by
multiplication by appropriate factors. The operators
$T^\epsilon _{bc}(I_{21})$ and $T^\epsilon _{bc}(I_{43})$ are given in
the basis $\{ | k,l\rangle ' \}$ by the same formulas (31).

We analytically continue formula (35) to the velues of $b$ and $c$ for
which the operators $T^\epsilon _{bc}(I_{32})$ give reducible representations
(that is, to those values of $b$ and $c$ which are excluded in Theorem 1). As
a result, we obtain the new operators ${\tilde T}^\epsilon _{bc}(I_{32})$
which gives with the operators
$T^\epsilon _{bc}(I_{21})$ and $T^\epsilon _{bc}(I_{43})$ an irreducible
representation ${\tilde T}^\epsilon _{bc}$ of $U'_q({\rm so}_{2,2})$. As in
the similar situation of the case of representations of Lie algebras (see
[14]), the reducible representations ${\tilde T}^\epsilon _{bc}$ and
$T^\epsilon _{bc}$ consist of the same irreducible components. It is more
convenient to look for irreducible components of the reducible representations
${\tilde T}^\epsilon _{bc}$. Depending on values of $b$ and $c$ we differ
5 cases.
\medskip

\noindent
{\sl Case 1:} Let $b$ be integral if $\epsilon =0$ and half-integral if
$\epsilon =1$, and let $c\ne c',c'+{\rm i}\pi /2h$, where $c'\in {\bf Z}$ if
$\epsilon =0$ and $c'$ is half-integral if $\epsilon =1$. (Recall that we
assumed that $b\ge 1/2$ and $0\le {\rm Im}\, c< \pi {\rm i}/h$.) Then the
second and third summands on the right hand side of (35) vanish for
$k=-b$, $k=b-1$ and $k=b$, $k=-b+1$, respectively. By the same reasoning as in
the case of representations of Lie algebras (see, for example, [14],
Chap. 7) we find that these vanishings of summands lead to appearence
of the irreducible invariant subspaces $H^0_{bc}$, $H^+_{bc}$ and $H^-_{bc}$
of the representation space $H_\epsilon$ spanned by all the basis vectors
$| k,l\rangle '$ with $-b<k<b$, $k\ge b$ and $k\le -b$, respectively.
(Note that if $b=1/2$, then the subspace $H^0_{bc}$ is empty.) We denote
the corresponding irreducible subrepresentations of
${\tilde T}^\epsilon _{bc}$ by $D^0_{bc}$, $D^+_{bc}$ and $D^-_{bc}$,
respectively. We have
$$
{\tilde T}^\epsilon _{bc}\sim D^0_{bc}\oplus D^+_{bc}\oplus D^-_{bc}\ \
{\rm if}\ \ b>1/2\ \ {\rm and}\ \
{\tilde T}^\epsilon _{1/2,c}\sim  D^+_{1/2,c}\oplus D^-_{1/2,c}.
$$
Note that if $b=1$ then the subspace $H^0_{bc}$ is spanned by the basis
vectors $| 0,l\rangle '\equiv | l\rangle$. In this case the operator
$D^0_{1,c}(I_{32})$ is given by the formula
$$
D^0_{1,c}(I_{32}) |l\rangle =([l]^2/[2l]^2)\{ ((q^{l+1}+q^{-l-1})
(q^{-l}+q^{l})[l+c][l-c+1])^{1/2}| l+1\rangle $$
$$
+((q^{l}+q^{-l})
(q^{-l+1}+q^{l-1})[l-c][l+c-1])^{1/2}| l-1\rangle \} . \eqno (36)
$$

\noindent
{\sl Case 2}: Let $c$ be integral if $\epsilon =0$ and half-integral if
$\epsilon =1$, and let $b\ne b',b'+{\rm i}\pi /2h$, where $b'\in {\bf Z}$ if
$\epsilon =0$ and $b'$ is half-integral if $\epsilon =1$. In this case
the representation space $H_\epsilon$ has the irreducible invariant subspaces
$H^0_{bc}$, $H^+_{bc}$ and $H^-_{bc}$ spanned by the basis vectors
$| k,l\rangle '$ with $-c<l<c$, $l\ge c$ and $l\le -c$, respectively.
We denote the corresponding irreducible subrepresentations of
${\tilde T}^\epsilon _{bc}$ by $F^0_{bc}$, $F^+_{bc}$ and $F^-_{bc}$,
respectively. We have
$$
{\tilde T}^\epsilon _{bc}\sim F^0_{bc}\oplus F^+_{bc}\oplus F^-_{bc}\ \
{\rm if}\ \ b>1/2\ \ {\rm and}\ \
{\tilde T}^\epsilon _{b,1/2}\sim  F^+_{b,1/2}\oplus F^-_{b,1/2}.
$$
The representation $F^0_{b,1}$ acts on the space spanned by the vectors
$| k,0\rangle \equiv |k\rangle$ and the operator $F^0_{b,1}(I_{32})$ is
given by the formula similar to the relation (36).
\medskip

\noindent
{\sl Case 3}: Let $b$ be as in Case 1 and let $c$ be of the form
$c= c'+{\rm i}\pi /2h$,
where $c'\in {\bf Z}$ if $\epsilon =0$ and $c'$ is half-integral if
$\epsilon =1$ and such that $c'>b$. Then the second
summand on the right hand side of (35) vanishes when $k=-b$, $k=b-1$,
$k=-c$, $k=c-1$ and the third summand vanishes when $k=b$, $k=-b+1$,
$k=c$, $k=-c+1$. This leads to appearence of five irreducible invariant
subspaces $H^0_{bc}$, $H^+_{bc}$, $H^{++}_{bc}$, $H^-_{bc}$, $H^{--}_{bc}$
in $H_\epsilon$ spanned by all the basis vectors $| k,l\rangle '$ with
$-b<k<b$, $b\le k<c'$, $k\ge c'$, $-c'<k\le -b$ and $k\le -c'$,
respectively. (Note that if $b=1/2$, then subspace $H^0_{bc}$ is empty.)
We denote the corresponding irreducible subrepresentations of
${\tilde T}^\epsilon _{bc}$  by
$Q^0_{bc}$, $Q^+_{bc}$, $Q^{++}_{bc}$, $Q^-_{bc}$, $Q^{--}_{bc}$.
We have
$$
{\tilde T}^\epsilon _{bc}\sim
Q^0_{bc}\oplus Q^+_{bc}\oplus Q^{++}_{bc}\oplus Q^-_{bc}\oplus Q^{--}_{bc}
\ \ \ {\rm if}\ \ \ b>1/2, $$
$$
{\tilde T}^\epsilon _{1/2,c}\sim Q^+_{1/2,c}\oplus Q^{++}_{1/2,c}\oplus
Q^-_{1/2,c}\oplus Q^{--}_{1/2,c}.
$$
If $b=1$, then the subspace $H^0_{bc}$ is spanned by the vectors
$|0,l\rangle '\equiv | l\rangle$ and the operator $Q^0_{bc}(I_{32})$ is
given by the formula similar to the relation (36). If $c'=b+1$, then the
subspace $H^+_{bc}$ is spanned by the basis vectors $|b,l\rangle '\equiv
|l\rangle$ and the operator $Q^+_{bc}(I_{32})$ is given by
$$
Q^+_{bc}(I_{32})| l\rangle =-([b+l][b-l])/([2(b+l)][2(b-l)])$$
$$
\times \{ ((q^{l+b}+q^{-l-b}) (q^{b-l-1}+q^{l-b+1})
(q^{l+c'}+q^{-l-c'}) (q^{c'-l+1}+q^{l-c'-1}))^{1/2} |l+1\rangle $$
$$
((q^{l+b-1}+q^{-l-b+1}) (q^{b-l}+q^{l-b})
(q^{l+c'-1}+q^{-l-c'+1}) (q^{c'-l}+q^{l-c'}))^{1/2} |l-1\rangle \} .
\eqno (37)
$$
At $c'=b+1$ the subspace $H^-_{bc}$ is spanned by the basis vectors
$| -b,l\rangle '\equiv |l\rangle$ and the operator $Q^-_{bc}(I_{32})$
is given by the formula similar to the relation (37).

If $b$ and $c$ are as above and $c'=b$, then the subspaces $H^+_{bc}$ and
$H^-_{bc}$ are empty and we have
$$
{\tilde T}^\epsilon _{bc}\sim
Q^0_{bc}\oplus Q^{++}_{bc}\oplus Q^{--}_{bc}
 \ \ {\rm if}\ \  b>1/2 \ \ {\rm and}\ \
{\tilde T}^\epsilon _{1/2,c}\sim Q^{++}_{1/2,c}\oplus Q^{--}_{1/2,c}.
$$

\noindent
{\sl Case 4}: Let $c$ be as in Case 2 and let $b$ be of the form
$b=b'+{\rm i}\pi /2h$,
where $b'\in {\bf Z}$ if $\epsilon =0$ and $b'$ is half-integral if
$\epsilon =1$ and such that $b'>c$. In this case the representation space
$H_\epsilon$ has the irreducible invariant subspaces
$H^0_{bc}$, $H^+_{bc}$, $H^{++}_{bc}$, $H^-_{bc}$, $H^{--}_{bc}$
spanned by all the basis vectors $| k,l\rangle '$ with
$-c<l<c$, $c\le l<b'$, $l\ge b'$, $-b'<l\le -c$ and $l\le -b'$,
respectively. If $c=1/2$, then the subspace $H^0_{bc}$ is empty.
We denote the corresponding irreducible subrepresentations of
${\tilde T}^\epsilon _{bc}$  by
$R^0_{bc}$, $R^+_{bc}$, $R^{++}_{bc}$, $R^-_{bc}$, $R^{--}_{bc}$.
We have
$$
{\tilde T}^\epsilon _{bc}\sim
R^0_{bc}\oplus R^+_{bc}\oplus R^{++}_{bc}\oplus R^-_{bc}\oplus R^{--}_{bc}
\ \ \ {\rm if}\ \ \ c>1/2, $$
$$
{\tilde T}^\epsilon _{b,1/2}\sim R^+_{b,1/2}\oplus R^{++}_{b,1/2}\oplus
R^-_{b,1/2}\oplus R^{--}_{b,1/2}.
$$
If $c=1$, then the subspace $H^0_{bc}$ is spanned by the vectors
$|k,0\rangle '\equiv | k\rangle$ and the operator $R^0_{bc}(I_{32})$ is
given by the formula similar to the relation (36). If $b'=c+1$, then the
subspace $H^+_{bc}$ (the subspace $H^-_{bc}$)
is spanned by the basis vectors $|k,c\rangle '\equiv |k\rangle$
(respectively, by $|k,-c\rangle '\equiv |k\rangle$)
and the operator $R^+_{bc}(I_{32})$ (respectively,
the operator $R^-_{bc}(I_{32})$) is given by the formula of the type (37).

If $b$ and $c$ are as above and $b'=c$, then the subspaces $H^+_{bc}$
and $H^-_{bc}$ are empty and we have
$$
{\tilde T}^\epsilon _{bc}\sim
R^0_{bc}\oplus R^{++}_{bc}\oplus R^{--}_{bc}
 \ \ {\rm if} \ \ c>1/2\ \ {\rm and}\ \
{\tilde T}^\epsilon _{b,1/2}\sim R^{++}_{b,1/2}\oplus R^{--}_{b,1/2}.
$$

\noindent {\sl Case 5:}
Let $b$ and $c$ be integral if $\epsilon =0$ and half-integral if $\epsilon
=1$. (Note that according to our convention we assume that $b\ge 1/2$ and
$c\ge 1/2$.) Then the first, second, third and fourth summands on the
right hand side of (35) vanish for the appropriate values of $k$ and $l$.
This leads to the decomposition of $H_\epsilon$ into the direct sum of the
irreducible invariant subspaces $H^{\epsilon _1,\epsilon _2}_{bc}$,
$\epsilon _1, \epsilon _2=0,+,-$, which are spanned by the basis vectors
$| k,l\rangle '$ with
$$
-b<k<b,\ \ -c<l<c \ \ \ {\rm for} \ \ \ E^{00}_{bc}; \ \ \ \ \
-b<k<b,\ \ l\ge c \ \ \ {\rm for} \ \ \ E^{0+}_{bc}; $$
$$
-b<k<b,\ \ l\le -c \ \ \ {\rm for} \ \ \ E^{0-}_{bc}; \ \ \ \ \
k\ge b,\ \ -c<l<c \ \ \ {\rm for} \ \ \ E^{+0}_{bc}; $$
$$
k\ge b,\ \ l\ge c \ \ \ {\rm for} \ \ \ E^{++}_{bc}; \ \ \ \ \
k\ge b,\ \ l\le -c \ \ \ {\rm for} \ \ \ E^{+-}_{bc}; $$
$$
k\le -b,\ \ -c<l<c \ \ \ {\rm for} \ \ \ E^{-0}_{bc}; \ \ \ \ \
k\le -b,\ \ l\ge c \ \ \ {\rm for} \ \ \ E^{-+}_{bc}; $$
$$
k\le -b,\ \ l\le -c \ \ \ {\rm for} \ \ \ E^{--}_{bc}.
$$
(Note that the subspaces $E^{0,\epsilon _2}_{bc}$ are absent if $c=1/2$ and
the subspaces $E^{\epsilon _1,0}_{bc}$ are absent if $b=1/2$.)
The corresponding subrepresentations of ${\tilde T}^\epsilon _{bc}$ are
denoted by $E^{\epsilon _1,\epsilon _2}_{bc}$, respectively.
We have
$$
{\tilde T}^\epsilon _{bc}=\sum _{\epsilon _1,\epsilon _2=0,+,-}
\oplus E^{\epsilon _1,\epsilon _2}_{bc}\ \ \ {\rm if} \ \ \
c\ne 1/2,\ \ b\ne 1/2, $$
$$
{\tilde T}^\epsilon _{1/2,c}=\sum _{\epsilon _1=+,-}
\sum _{\epsilon _2=0,+,-}
\oplus E^{\epsilon _1,\epsilon _2}_{1/2,c}\ \ \ {\rm if} \ \ \
c\ne 1/2, $$
$$
{\tilde T}^\epsilon _{b,1/2}=\sum _{\epsilon _1=0,+,-}
\sum _{\epsilon _2=+,-}
\oplus E^{\epsilon _1,\epsilon _2}_{b,1/2}\ \ \ {\rm if} \ \ \
b\ne 1/2, $$
$$
{\tilde T}^\epsilon _{1/2,1/2}=E^{++} _{1/2,1/2}\oplus
E^{+-} _{1/2,1/2}\oplus E^{-+} _{1/2,1/2}\oplus E^{--} _{1/2,1/2}.
$$
Clearly, the irreducible representations $E^{00}_{bc}$ and only they
are finite dimensional.
\medskip

\noindent {\bf Theorem 2.}
(1) {\it The irreducible representations $T^\epsilon _{bc}$ and irreducible
components of reducible representations $T^\epsilon _{bc}$ lead to the
following classes of irreducible representations of $U'_q({\rm so}_{2,2})$:}
\medskip

(a) {\it The representations of Theorem 1};

(b) {\it The representations $D^{0}_{bc}$, $D^{+}_{bc}$, $D^{-}_{bc}$,
where $b\in {\frac 12}{\bf Z}$, $b\ge 1/2$};

(c) {\it The representations $F^{0}_{bc}$, $F^{+}_{bc}$, $F^{-}_{bc}$,
where $c\in {\frac 12}{\bf Z}$, $c\ge 1/2$};

(d) {\it The representations $Q^{0}_{bc}$, $Q^{+}_{bc}$, $Q^{++}_{bc}$,
$Q^{-}_{bc}$, $Q^{--}_{bc}$, where $b\in {\frac 12}{\bf Z}$ and
$c=c'+{\rm i}\pi /2h$, $c'\in {\frac 12}{\bf Z}$, $b,c'\ge 1/2$};

(e) {\it The representations $R^{0}_{bc}$, $R^{+}_{bc}$, $R^{++}_{bc}$,
$R^{-}_{bc}$, $R^{--}_{bc}$, where $c\in {\frac 12}{\bf Z}$ and
$b=b'+{\rm i}\pi /2h$, $b'\in {\frac 12}{\bf Z}$, $b',c\ge 1/2$};

(f) {\it The representations $E^{\epsilon _1,\epsilon _2}_{bc}$,
$\epsilon _1,\epsilon _2=0,+,-$, where
$b,c\in {\frac 12}{\bf Z}$, $b,c\ge 1/2$}.
\smallskip

\noindent
{\it Every irreducible representation of $U'_q({\rm so}_{2,2})$, which is
equivalent to some irreducible representation $T^\epsilon _{bc}$ or to an
irreducible component of some reducible representation $T^\epsilon _{bc}$,
is equivalent to one of the representations of classes (a)-(f)}.
\smallskip

(2) {\it Between representations of classes (a)-(f) there exist no
equivalence relations except for relations which are compositions of
the relations (32)-(34).}
\medskip

\noindent {\sl Proof.}
Proof of the assertion (1) is given above. The assertion (2) for the
representations of class (a) is given above. Absence of other equivalence
relations follows from the fact that representations of any other pair
from (a)-(f) have non-coinciding spectra of the operators $T(I_{21}$
and $T(I_{43})$.
\bigskip

\centerline{{\sc 8. Irreducible $*$-representations of} $U'_q({\rm so}_{2,2})$}
\medskip

The aim of this section is to give the classification of $*$-representations
in the set of irreducible representations of $U'_q({\rm so}_{2,2})$
from Theorem 2. This classification is derived by the calculations described,
for example, in [15], Sect. 6.4.6. For this reason, we only formulate the
final result.
\medskip

\noindent {\bf Theorem 3.} {\it The following representations in the set of
irreducible representations of Theorem 2 are $*$-representations of
$U'_q({\rm so}_{2,2})$:}
\medskip

(1) {\it The representations $T^\epsilon _{bc}$ when $b={\rm i}\rho +1/2$,
$c={\rm i}\rho '+1/2$, $\rho ,\rho '\in {\bf R}$;}
\smallskip

(2) {\it The representations $T^\epsilon _{bc}$, $\epsilon =0$, when $b$ is
in one of the intervals $(0,1/2]$, $(0+{\rm i}\pi/2h, 1/2+{\rm i}\pi/2h]$
and $c$ is in one of these integrals;}
\smallskip

(3) {\it The representations $T^\epsilon _{bc}$, $\epsilon =0$, when
$b={\rm i}\rho +1/2$, $\rho \in {\bf R}$, and $c$ is
in one of the intervals $(0,1/2]$, $(0+{\rm i}\pi/2h, 1/2+{\rm i}\pi/2h]$
or when $c={\rm i}\rho +1/2$, $\rho \in {\bf R}$, and $b$ is in one of
these integvals;}
\smallskip

(4) {\it The representations $D^+_{bc}$, $D^-_{bc}$, $D^0_{1,c}$ when
$c={\rm i}\rho +1/2$, $\rho \in {\bf R}$, or when $c$ belongs
to one of the intervals $(0,1/2]$, $(0+{\rm i}\pi/2h, 1/2+{\rm i}\pi/2h]$;}
\smallskip

(5) {\it The representations $F^+_{bc}$, $F^-_{bc}$, $F^0_{b,1}$ when
$b={\rm i}\rho +1/2$, $\rho \in {\bf R}$, or when $b$ belongs
to one of the intervals $(0,1/2]$, $(0+{\rm i}\pi/2h, 1/2+{\rm i}\pi/2h]$;}
\smallskip

(6) {\it The representations $Q^{++}_{bc}$, $Q^{--}_{bc}$,
$R^{++}_{bc}$, $R^{--}_{bc}$;}
\smallskip

(7) {\it The representations $Q^{+}_{bc}$, $Q^{-}_{bc}$ with
$c=b+1+{\rm i}\pi /2h$ and the representations
$R^{+}_{bc}$, $R^{-}_{bc}$ with $b=c+1+{\rm i}\pi /2h$;}
\smallskip

(8) {\it The representations $Q^0_{bc}$ with
$c=b+{\rm i}\pi /2h$ and the representations
$R^0_{bc}$ with $b=c+{\rm i}\pi /2h$;}
\smallskip

(9) {\it The representations $E^{\epsilon _1,\epsilon _2}_{bc}$,
$\epsilon _1,\epsilon _2=+,-$, and the representations $E^{0,\pm}_{1,c}$,
$E^{\pm ,0} _{b,1}$.}
\medskip

The representations of class (1) are called $*$-representations of the
{\it principal series}.
The representations of class (2) are called $*$-representations of the
{\it supplementary series}.
The representations of class (3) belong to the mixed (principal--supplementary)
series.
The representations of $D^{\pm}_{bc}$, $F^{\pm}_{bc}$ of classes (4) and (5)
also belong to the mixed (discrete--principal and discrete--supplementary)
series.
The representations of class (6) are called $*$-representations of the
{\it discrete series}.
The representations of $D^0_{1,c}$ and $F^0_{b,1}$ from
classes (4) and (5),
the representations of classes (7) and (8) and
the representations $E^{0,\pm}_{1,c}$, $E^{\pm ,0}_{b,1}$ are called
{\it ladded} $*$-representations. The operators $T(I_{32})$ for these
representations are of the types (36) and (37).
\bigskip

\centerline{\sc 9. Conclusions}
\medskip

{\bf 1.} In Theorems 2 and 3 we described sets of irreducible
representations and irreducible $*$-representations of the algebra
$U'_q({\rm so}_{2,2})$, respectively. The unsolved problem is the
following: Do these sets of representations exhaust (up to equivalence)
all irreducible representations and $*$-representations of this algebra?
\smallskip

{\bf 2.} The algebra $U'_q({\rm so}_{2,2})$ is a $q$-deformation of the
universal enveloping algebra $U({\rm so}_{2,2})$ of the Lie algebra
so$_{2,2}$. It is well-known that so$_{2,2}$ is the direct sum of two Lie
algebras so$_{2,1}$: ${\rm so}_{2,2}={\rm so}_{2,1}\oplus {\rm so}_{2,1}$.
Therefore, irreducible representations and irreducible $*$-representations
of so$_{2,2}$ are easily determined by the representations and
$*$-representations of the Lie algebra so$_{2,1}$: any
irreducible representation ($*$-representation) $T$ of so$_{2,2}$ is a
direct product of two such representations of so$_{2,1}$. This trivially
gives the classification of
irreducible representations and irreducible $*$-representations of
so$_{2,2}$. Comparing these classifications with the irreducible
representations of $U'_q({\rm so}_{2,2})$ we see that the representations
of classes (d) and (c) have no analogue for the Lie algebra so$_{2,2}$.
Similarly, many $*$-representations of Theorem 3 (for example, the
representations with $b$ or $c$ lying in the interval $(0+{\rm i}\pi /2h,
1/2+{\rm i}\pi /2h]$ have no analogue for so$_{2,2}$.

{\bf 3.} Comparing the classification of irreducible $*$-representations
of the $q$-deformed algebras $U'_q({\rm so}_{2,1})$ and
$U'_q({\rm so}_{3,1})$ (see [2]) with irreducible $*$-representations
of $U'_q({\rm so}_{2,2})$ in Theorem 3 we make the following conclusions:
\smallskip

(a) The algebra $U'_q({\rm so}_{2,2})$ (unlike the algebras
$U'_q({\rm so}_{2,1})$ and $U'_q({\rm so}_{3,1})$) has no strange series of
irreducible $*$-representations.

(b) The algebras $U'_q({\rm so}_{2,1})$ and $U'_q({\rm so}_{3,1})$
(unlike the algebra $U'_q({\rm so}_{2,2})$) have no mixed series
of irreducible $*$-representations.
\smallskip

These conclusions say that under the transition from the $q$-deformed
algebras corresponding to Lie algebras of rank 1 to the
$q$-deformed algebras corresponding to Lie algebras of higher ranks
we obtain a quantitative difference in their representation theories. It
is not known now if this property is valid for representations of the
Drinfeld-Jimbo algebras.
\bigskip

{\sc Acknowledgement.} This research was supported
in part by CRDF Grant UP1-309 and by DFFD Grant 1.4/206.

\centerline{\sc Bibliography}
\medskip

1. A. M. Gavrilik and A. U. Klimyk, {\it $q$-Deformed orthogonal and
pseudo-orthogonal algebras and their representations},
Lett. Math. Phys. {\bf 21} (1991), 215--220.

2. A. M. Gavrilik and A. U. Klimyk, {\it Representations of $q$-deformed
algebras $U_q({\rm so}_{2,1})$ and $U_q({\rm so}_{3,1})$},
J. Math. Phys. {\bf 35} (1994), 3670--3686.

3. M. Noumi, {\it Macdonald's symmetric polynomials as zonal spherical
functions on quantum homogeneous spaces}, Adv. in Math. {\bf 123} (1996),
16--77.

4. Yu. S. Samoilenko and L. B. Turowska, {\it Semilinear relations and
$*$-representations of deformations of $SO(3)$}, Reps. Math. Phys., in press.

5. O. V. Bagro, S. A. Kruglyak, {\it Representations of Fairlie algebra},
preprint, Kiev, 1996.

6. A. M. Gavrilik, {\it Representations of $U_q({\rm so}_{4})$ and
$U_q({\rm so}_{3,1})$}, Teoret. Matem. Fiz. {\bf 95} (1993), 251--257.

7. V. G. Drinfeld, {\it Hopf algebras and the quantum Yang-Baxter
equation}, Soviet Math. Dokl. {\bf 32} (1985), 667--671.

8. M. Jimbo, {\it A $q$-analogue of $U(gl(N+1))$ and the quantum
Yang-Baxter equation}, Lett. Math. Phys. {\bf 10} (1985), 63--69.

9. A. U. Klimyk, K. Schm\"udgen,  {\it Quantum Groups and Their
Representations}, Springer, Berlin, 1998.

10. I. M. Gel'fand, M. L. Tsetlin, {\it Finite dimensional
representations of the group of orthogonal matrices},
Dokl. Akad. Nauk SSSR {\bf 71} (1950), 1017--1020.

11. D. B. Fairlie, {\it Quantum deformations of $SU(2)$},
J. Phys. A: Math. Gen. {\bf 23} (1990), L183--L187.

12. I. I. Kachurik, A. U. Klimyk, {\it Representations of the $q$-deformed
algebra $U'_q({\rm so}_4)$}, J. Phys. A {\bf 27} (1994), 7087--7097.

13. G. Gasper G, M. Rahman,  {\it Basic Hypergeometric Functions},
Cambridge Univ. Press, Cambridge, 1991.

14. A. U. Klimyk,  {\it Matrix Elements and Clebsch-Gordan Coefficients of
Group Representations}, Naukova Dumka, Kiev, 1979.

15. N. Ja. Vilenkin, A. U. Klimyk,  {\it Representation of Lie Groups and
Special Functions}, Kluwer, Dordrecht, 1991.
\bigskip

\end{document}